\documentclass[10 pt,leqno]{amsart}

\usepackage {amsfonts}
\usepackage{amsthm}
\usepackage{amssymb}
\usepackage{latexsym}
\usepackage{amsmath}
\usepackage{mathrsfs}
\usepackage{pifont}

\newcommand{\la} {\langle}
\newcommand{\ra} {\rangle}
\newcommand{\lla} {\left\langle}
\newcommand{\rra} {\right\rangle}

\newcommand{\QSproj} {\Delta^{QS}}

\newcommand{\bc} {\mathbb C}
\newcommand{\be} {\mathbb E}
\newcommand{\bn}{\mathbb N}
\newcommand{\br}{\mathbb R}
\newcommand{\Rplus}{\mathbb R_+}

\newcommand{\alg} {\mathsf A}

\newcommand{\Ran}{\rm{Ran}}
\newcommand{\Lin}{\rm{Lin}}
\newcommand{\clLin}{\overline{\rm{Lin}}}

\newcommand {\Op} {\mathcal{O}}

\newcommand{\Hil}{\mathsf{H}}
\newcommand{\hil}{\mathsf{h}}
\newcommand{\Kil}{\mathsf{K}}
\newcommand{\kil}{\mathsf{k}}
\newcommand{\Vil}{\mathsf{V}}

\newcommand{\noise}{\mathsf{k}}

\newcommand{\Wil}{\mathsf{W}}

\newcommand{\wh}{\widehat}

\newcommand{\Com}{\Delta}
\newcommand{\Cou}{\epsilon}

\newcommand{\ve}{\varepsilon}
\newcommand{\Fock}{\mathcal F}
\newcommand{\FFock}{\mathcal F}
\newcommand{\fhat}{\hat{f}}
\newcommand{\ghat}{\hat{g}}

\newcommand{\khat}{\widehat{\noise}}
\newcommand{\kilhat}{\widehat{\kil}}

\newcommand{\Exps}{\mathcal{E}}
\newcommand{\Step}{\mathbb{S}}

\theoremstyle{plain}
\newtheorem{propn}{Proposition}[section]
\newtheorem{thm}[propn]{Theorem}
\newtheorem{lemma}[propn]{Lemma}

\theoremstyle{definition}
\newtheorem{defn}[propn]{Definition}

\newtheorem*{remarks}{Remarks}

\newcommand{\Pto}{\!\to\!}

\newcommand{\ProchilE}{\mathbb{P}(\hil\odot\Exps)}

\newcommand{\ProckilhatE}{\mathbb{P}(\kilhat\odot\Exps)}
\newcommand{\ProcE}{\mathbb{P}(\Exps)}
\newcommand{\ProcEk}{\mathbb{P}(\Exps_{\kil})}
\newcommand{\ProcVtoW}{\mathbb{P}\big(\Vil; \Wil,\Exps\big)}

\newcommand{\ProcVtoV}{\mathbb{P}\big(\Vil;\Vil ,\Exps\big)}

\newcommand{\ProcVtoC}{\mathbb{P}\big(\Vil ;\Exps\big)}

\newcommand{\ProcAtoC}{\mathbb{P}\big(\alg;\Exps\big)}
\newcommand{\PProcAtoCk}{\mathbb{P}\big(\alg;\Exps_{\kil}\big)}
\newcommand{\AtoCProc}{\mathbb{P}\big(\alg;\Exps_{\kil_0}\big)}
\newcommand{\CdProcAtoCd}{\mathbb{P}\big(\alg;\Exps_{\kil}\big)}

\newcommand{\AwtProcwtAtoC}{\mathbb{P}\big(\wt{\alg};\Exps\big)}

\newcommand{\iotProcAtoCiota}{\mathbb{P}\big(\alg;\Exps^\iota\big)}

\newcommand{\salvation}{\mathbb{P}\big(\alg;\Exps_{\kil_0}\big)}

\newcommand{\ot}{\otimes}

\newcommand{\ida}{1_{\alg}}
\newcommand{\id}{\text{id}}

\newcommand{\idB}{\textrm{id}_{B(\khat)}}
\newcommand{\idf}{I_{\Fock}}
\newcommand{\wt}{\widetilde}

\numberwithin{equation}{section}

\begin{document}
\author{Adam G. Skalski}
\footnote{\emph{Permanent address of the author}.
Department of Mathematics, University of \L\'{o}d\'{z}, ul.
Banacha 22, 90-238 \L\'{o}d\'{z}, Poland.}
\title{
%Dilation of quantum stochastic \\ convolution cocycles on a $C^*$-bialgebra}
Completely positive quantum stochastic convolution cocycles and their dilations}
\address{Department of Mathematics and Statistics,  Lancaster University,
Lancaster, LA1 4YF}
\email{a.skalski@lancaster.ac.uk}


\begin{abstract}
Stochastic generators of completely positive and contractive quantum
stochastic convolution cocycles on a $C^*$-hyperbialgebra are
characterised. The characterisation is used to obtain dilations and
stochastic forms of Stinespring decomposition for completely positive
convolution cocycles on a $C^*$-bialgebra.
\end{abstract}

\maketitle

Stochastic (or Markovian) cocycles on operator algebras are basic objects
of interest in quantum probability (\cite{Acc}) and have been extensively
investigated using quantum stochastic analysis (see~\cite{lect}). There is
also a well-developed theory of quantum L\'evy processes, that is
stationary, independent-increment, *-homomorphic processes on a
*-bialgebra (see~\cite{schu},~\cite{franz} and references therein).
Close examination of these two directions has naturally led to the
notion of \emph{quantum stochastic convolution cocycle} on a quantum group
(or, more generally, on a coalgebra), as introduced and investigated in
\cite{LSqscc1} in an algebraic context, and in~\cite{LSqscc2} in the
analytic context of compact quantum groups. The main results have been
summarised in~\cite{LSbedlewo}.
Recent years have also seen an increased interest in the noncommutative
generalisation of classical hypergroups (\cite{hyp}), initiated by
Chapovsky and Vainerman (\cite{ChaV}) and continued, for example, in
the papers~\cite{Kal1} and~\cite{Kal2}. Compact quantum hypergroups differ
from compact quantum groups in that their coproduct need not be
multiplicative. However, it remains completely positive, which
makes compact quantum hypergroups, or more generally
$C^*$-hyperbialgebras, an appropriate category for the consideration of
completely positive quantum stochastic convolution cocycles in a
topological context (for the purely algebraic case see~\cite{UweSch}).
These cocycles may be viewed as natural counterparts of stationary,
independent-increment processes on hypergroups. In \cite{LSqscc2} it is
shown that, under certain regularity conditions, they satisfy
coalgebraic quantum stochastic differential equations.

The aim of this paper is to
% stochastic generators of completely positive and contractive quantum stochastic convolution
%cocycles on a $C^*$-hyperbialgebra, and to use it to \end{comment}
prove dilation theorems for quantum stochastic convolution cocycles
on a $C^*$-bialgebra. To this end it is first necessary to establish
the detailed structure of the stochastic generators of completely positive and contractive
convolution cocycles. We give a direct derivation of this exploiting ideas used in the analysis of 
 standard quantum stochastic cocycles with finite-dimensional noise space (\cite{lp}).
Once the structure of generators is known,
% detailed description of the structure of stochastic generators is established,
one may consider question of dilating completely positive convolution cocycles to
$^*$-homomorphic ones. In the context of \emph{standard} quantum stochastic cocycles
this problem was treated in \cite{dilate} and \cite{Stine} (see also \cite{Slava}).
In the first of these papers it was shown that every Markov-regular completely positive
and contractive cocycle arises
as the image of a $^*$-homomorphic cocycle under a vacuum conditional expectation which
averages out some dimensions of the quantum noise.
In the second every Markov-regular completely positive and contractive cocycle was shown
to be realisable as a composition of a $^*$-homomorphic cocycle with conjugation by
a contraction operator process. This may be seen as a stochastic Stinespring decomposition.
 In this paper using the techniques of Goswami, Lindsay, Sinha and Wills we obtain  analogous
results for convolution cocycles on $C^*$-bialgebras. Multiplicativity of the coproduct is
necessary here for obtaining dilations to $^*$-homomorphic cocycles.

An alternative approach to the one presented in this paper
would be to exploit more directly theorems known for standard quantum stochastic
convolution cocycles and properties of the $R$-map introduced in \cite{LSqscc2}.
%The structure of stochastic generators of standard quantum
%stochastic cocycles, which are completely positive, contractive and
%Markov-regular, has been considered in a number of papers
%(\cite{Slava},~\cite{lp},~\cite{lw1},~\cite{lwblms}) and is now
%well-understood.
In that paper the general form of the 
stochastic generators of completely positive and contractive convolution cocycles
was determined by using a particular
representation of the $C^*$-bialgebra in question and appealing directly to the results of 
\cite{lp},~\cite{lw1} and ~\cite{lwblms}; similar methods may be further used to obtain the dilation results presented here.
One drawback of such an approach is that it involves using the deep Christensen-Evans theorem on quasi-innerness of derivations on represented $C^*$-algebras. Another is the necessity to reformulate reformulation of the results of \cite{dilate}
and \cite{Stine} in coordinate-free quantum stochastic calculus. This is also necessary for overcoming separability assumptions on the noise dimension spaces. Finally the von Neumann
algebraic framework used in \cite{Stine} would require nontrivial modifications.
In sum the method presented here has the advantage of being more elementary.

%Dilation of such cocycles to *-homomorphic quantum
%stochastic cocycles was achieved in~\cite{dilate} and~\cite{Stine}.
%In the first of these papers it was shown that every Markov-regular,
%completely positive, contractive cocycle arises as the image of a
%*-homomorphic cocycle under a vacuum conditional expectation which
%averages out some dimensions of the quantum noise.
%5In the second every Markov-regular, completely positive, contractive
%cocycle was shown to be realisable as a composition of a *-homomorphic
%cocycle with conjugation by a contraction process governed by a quantum
%stochastic differential equation whose coefficients evolve in time
%according to the homomorphic process.

%Here completely positive, contractive, Markov-regular quantum stochastic
%onvolution cocycles on a $C^*$-hyperbialgebra are considered. A direct
%and elementary proof of the detailed structure of their stochastic
%generators is given, following ideas used in the original paper
%characterising the generators of the corresponding class
%of standard cocycles (\cite{lp}). This approach avoids appeal to the
%Christensen-Evans theorem on quasi-innerness of derivations on represented
%$C^*$-algebras. Using techniques from the papers of Goswami, Lindsay,
%Sinha and Wills, analogous dilation results and stochastic Stinespring
%decompositions for convolution cocycles on $C^*$-bialgebras are obtained.
%

The plan of the paper is as follows.
The first section contains the notation and elements of quantum stochastic
analysis and operator space theory needed.
 In the second section $C^*$-hyperbialgebras are
defined and the well-known technique of obtaining them from
$C^*$-bialgebras via a noncommutative conditional expectation is recalled.
Basic facts concerning quantum stochastic convolution
cocycles and the structure of their stochastic generators in the completely positive and
$^*$-homomorphic cases are also included here.
In the third section a more detailed description of the
stochastic generators of Markov-regular, completely positive, contractive convolution
cocycles, in terms of a certain tuple of objects, is derived.
Dilations of such convolution cocycles on a
$C^*$-bialgebra to *-homomorphic convolution cocycles are given in the
fourth section, and the fifth section contains a stochastic Stinespring decomposition.

\section{Preliminaries}

In this section we introduce our notations and review results from quantum
stochastic analysis relevant to the rest of the paper. We shall usually
abbreviate
\emph{quantum stochastic} to QS,
\emph{completely positive} to CP and
\emph{completely positive, contractive} to CPC.

\subsection*{Notation}
All vector spaces in this paper are complex and inner products are linear
in their second argument. Algebraic tensor products are denoted by
$\odot$.

Let $\hil$ be a Hilbert space.
Ampliations are denoted
\begin{equation*}
%\label{iota hil}
\iota_{\hil}: B(\Hil)\to B(\Hil\ot\hil), \ T\mapsto T\ot I_{\hil},
\end{equation*}
and each vector $\xi\in\hil$ defines operators
\begin{equation*}
%\label{E-notation}
E_\xi : \Hil \to \Hil\ot\hil, \ v\mapsto v\ot\xi \text{ and }
E^{\xi} = (E_\xi)^*,
\end{equation*}
generalising Dirac's bra-ket notation:
\begin{equation*}
%\label{bra-ket}
E_\xi = I_{\Hil} \ot |\xi\ra \text{ and } E^\xi = I_{\Hil} \ot \la\xi |.
\end{equation*}
The particular Hilbert space $\Hil$ will always be clear from
the context.
For a subspace $E$ of $\hil$, $\Op (E)$ will denote the vector space of
linear operators in $\hil$ with domain $E$.

Finally, for a function $f:\br_+ \to \hil$ and subinterval $I$ of $\br_+$,
$f_I$ denotes the function $\br_+\to\hil$ which agrees with $f$ on $I$ and
is zero outside $I$ (cf.\ standard indicator-function notation).This
convention also applies to vectors, by viewing them as constant
functions --- for example
\[
\xi_{[s,t[}, \text{ for } \xi\in\hil  \text{ and } 0\leq s<t.
\]

\subsection*{Matrix spaces}
For an introduction to the theory of operator spaces we refer to
\cite{ERuan}. For this paper it is sufficient to work with concrete
operator spaces, that is closed subspaces of the space $B(\hil;\hil')$
of all bounded linear operators acting between Hilbert spaces
$\hil$ and $\hil'$.
The spatial/minimal tensor product of operator spaces is denoted by
$\otimes$, and when $\Vil, \Wil$ are operator spaces $CB(\Vil;\Wil)$
denotes the space of all completely bounded maps
from $\Vil$ to $\Wil$.

We need the concept of matrix spaces introduced by Lindsay and Wills in
\cite{lwblms}.
Let $\Vil\subset B(\Kil)$ be an operator space and let $\hil$ be a
supplementary Hilbert space. The operator space:
\[
M_{\hil} (\Vil) =
\{ T \in B(\Kil \ot \hil): E^{\xi'}TE_\xi \in \Vil
\text{ for all } \xi,\xi' \in \hil\}
\]
is called the \emph{$\hil$-matrix space} over $\Vil$. It is easy to see
that $M_{\hil} (\Vil)$ contains the spatial tensor product $\Vil
\ot B(\hil)$. Whenever $\Wil \subset B(\Hil)$ is another operator space,
and $\phi \in CB(\Vil;\Wil)$, the map
$\phi \ot \text{id}_{B(\hil)}$ extends uniquely to a completely bounded map
$\phi^{(\hil)}:  M_{\hil} (\Vil) \to   M_{\hil} (\Wil)$ satisfying
\[
E^{\xi'} \big(\phi^{(\hil)}(T)\big) E_\xi = \phi (E^{\xi'} T E_\xi),
\]
for all $T \in  M_{\hil} (\Vil)$, $\xi,\xi' \in \hil$. The map
$\phi^{(\hil)}$
is called the \emph{$\hil$-lifting} of $\phi$.

%NOT USED
%The \emph{$\hil$-column space} and  \emph{$\hil$-row space} over $\Vil$,
%are defined analogously, and are denoted respectively by
%$C_{\hil}(\Vil)\subset B(\Kil; \Kil \ot \hil)$ and $R_{\hil} (\Vil)$.
%Completely bounded maps lift to those as well.

\subsection*{Fock space notations and QS processes}

Let $\kil$ be a Hilbert space, called the \emph{noise dimension space}.
Then $\Fock_{\kil}$ denotes
the symmetric Fock space over $L^2(\br_+;\kil)$.
Exponential vectors in $\Fock_{\kil}$ are written $\ve(f)$, $f\in
L^2(\br_+;\kil)$.
The CCR flow of index $\kil$, defined in terms of
the second quantised shift on $L^2(\br_+;\kil)$,
is denoted $\sigma =\big(\sigma_t\big)_{t \geq 0}$.
Define
\[
\Step_{\kil} := \text{Lin}\{d_{[0,s[}: d \in \kil, s \in \br_+\}
\]
and a corresponding subspace of $\Fock$:
\[
\Exps_{\kil} := \text{Lin} \{\ve(f): f \in \Step_{\kil}\}.
\]
When the space $\kil$ is clear from the context we will simply write $\FFock$,
$\Step$ and $\Exps$.
Elements of $\Exps$ will play the role of test functions.
For a subspace $E$ of $\kil$ the following notation will be employed:
\begin{equation*}
\wh{E} := \Lin\big\{\wh{c}: c\in E\big\},\text{ where }
\wh{c}:= \binom{1}{c}\in\wh{\hil}:=\bc \oplus \hil.
\end{equation*}
Two further useful notations are
\begin{equation} \label{QSproj}
e_0 = \binom{1}{0}\in\kilhat \text{ and }
\QSproj = \begin{bmatrix} 0& \\&  I_{\kil} \end{bmatrix} \in B(\kilhat).
\end{equation}

Let $\hil$ be an additional Hilbert space.
By an \emph{$\hil$-operator process} we understand a family
$X = \big(X_t\big)_{t \geq 0}$ of operators on $\hil \ot\Fock$,
each having the (dense) domain $\hil \odot \Exps$, being weak-operator
measurable in $t$ and adapted to the natural Fock-space
operator filtration.
Thus $X:\br_+ \to \Op (\hil \odot \Exps)$, $t\mapsto X_t\xi$ is weakly
measurable
for all $\xi\in\hil \odot \Exps$ and, for each $t\geq 0$, $\zeta \in \hil$,
$X_t(\zeta \ot \ve(f))= X(t)\big(\zeta \ot \ve(f|_{[0,t[})\big)\ot\ve(f|_{[t,\infty[})$
for some operator $X(t)\in\Op\big(\hil \odot \Exps_{[0,t[}\big)$, where
$\Exps_{[0,t[}$ is defined as $\Exps$ is, except that $\br_+$ is
replaced by $[0,t[$.
The linear space of all $\hil$-operator processes is denoted
 $\ProchilE$, or  $\ProcE$ if $\hil=\bc$, with subscripts on the
$\Exps$ when necessary for avoiding ambiguity.
A process $X\in \ProchilE$
is called \emph{weakly regular} if for each
$\xi,\xi' \in \hil\odot\Exps$,
the scalar-valued function
\[ t \mapsto \la \xi', X_t \xi \ra, \;\; t\in \br_+,\]
is locally bounded.
It is called \emph{bounded} if $X_t$ is a bounded operator for each $t\geq 0$
(in such a case $X_t$ is usually identified with its continuous extension
to the whole of $\hil\ot \Fock$).

Now let $\Vil$ and $\Wil$ be operator spaces with $\Wil \subset B(\hil)$.
A linear map $k$ from $\Vil$ to $\ProchilE$ is called a
\emph{process on} $\Vil$ \emph{with values in} $\Wil$ if, for each
$v\in\Vil$, $k(v)$ is a $\hil$-operator process and, for each $f,g\in
\Step$, $t \geq 0$ and $v \in\Vil$, the operator
$E^{\ve(f)} k_t(v) E_{\ve(g)}$ belongs to $\Wil$.
Here $\Wil$ will usually be either $\Vil$ or $\bc$.

The vector space of all processes on $\Vil$ with values in $\Wil$ is written
$\ProcVtoW$ (this corresponds to the notation
$\mathbb{P}\big(\Vil\Pto\Wil :\hil\odot\Exps\big)$ used in \cite{LSqscc2});
when $\Wil=\bc$ we simply write  $\ProcVtoC$.
We say that $k \in \ProcVtoW$ is
\emph{pointwise weakly regular} if each $k(v)$ ($v\in \Vil$)  is weakly
regular. It is \emph{completely bounded} if,
for each $v \in \Vil$, the process $k(v)$ is bounded and, for each $t\geq 0$,
the map $k_t:\Vil \to B(\hil\ot \FFock)$ is completely bounded.

\subsection*{QS differential equations and standard QS cocycles}

Let $\Vil,\Wil$ be operator spaces with $\Wil \subset B(\hil)$ for some
Hilbert space $\hil$. For maps $\theta\in CB(\Vil;\Wil)$
and  $\phi \in CB\big(\Vil; M_{\kilhat}(\Vil)\big)$
consider the {\it quantum stochastic differential equation}
\begin{equation}
\label{QSDE}
dk_t = \widehat{k}_t \circ \phi\, d\Lambda_t, \;\;\;
k_0 = \iota_{\Fock} \circ \theta.
\end{equation}
 By a \emph{weak solution} of this equation we understand a process
$k \in \ProcVtoW$ such that
\[
\big\langle \xi \ot \ve(f),
\big(k_t(v) - \theta(v) 1_{\Fock}\big) \eta \ot\ve (g)\big\rangle =
\int_0^t {\big\langle \xi \ot \ve(f),
k_s \big(E^{\wh{f}(s)}\phi(v)E_{\ghat(s)}\big)\eta\ot\ve(g)\big\rangle}ds
\]
for all $t \geq 0$, $v\in \Vil$, $\xi,\eta \in \hil$ and $f,g\in \Step$.
If there is a \emph{quantum stochastically integrable}
$\hil\ot\khat$-process $K$ on $\Vil$ (see \cite{lect}), with
domain $\hil\ot\kilhat\odot\Exps$, satisfying
\[
E^{\zeta'}K_t(v)E_{\zeta} =
k_t\big(E^{\zeta'}\phi(v)E_{\zeta}\big)
\]
for all $t\geq 0$, $\zeta,\zeta'\in \kilhat$ and  $v\in \Vil$,
then $k$ is called a \emph{strong solution}. The equation (\ref{QSDE}) has
a unique weakly regular weak solution, which is also a strong solution
(\cite{lw1}, \cite{LSqsde}). We denote it by $k^{\theta,\phi}$, or simply
$k^{\phi}$ if $\Vil = \Wil$ and $\theta=\text{id}_{\Wil}$.

%\begin{defn}
A completely bounded process $k\in \ProcVtoV$ is called a \emph{standard QS cocycle on $\Vil$}
if, for $s,t \geq 0$,
\begin{equation*}
\label{coc}
 k_{t+s} = \hat{k}_t \circ \wt{\sigma}_s \circ k_s\;\;  \text{and}\;\;
 k_0 = \iota_{\FFock} \circ \text{id}_{\Vil}
\end{equation*}
where $\hat{k}_t$ denotes an $\Fock_{[s, \infty[}$-lifting of $k_t$ and $\wt{\sigma}_s = \text{id}_{\hil} \ot
\sigma_s$.
It is said to be \emph{Markov-regular}
if its \emph{Markov semigroup} $P: \Rplus\to B(\Vil)$, defined by
\[ P_t(v) = E^{\ve(0)} k_t(v) E_{\ve(0)}, \;\;\; t\geq 0, v \in \Vil,\]
is norm-continuous.
%\end{defn}
Whenever $\phi\in CB(\Vil; M_{\kilhat}(\Vil))$, the process $k^{\phi}$ is
a Markov-regular \emph{weak standard QS cocycle} (see \cite{lwjfa},
\cite{lect}; note however that what is here called a weak standard QS cocycle,
there is simply called a quantum stochastic cocycle).

%
%The following theorem originates in the earlier work of Lindsay and
%Parthasarathy (\cite{lp}) --- for the form of the generator
%see also \cite{Slava}.
%
%
%\begin{thm} [\cite{lw1}, \cite{lwblms}] \label{CPstandard}
%Let $\alg\subset B(\hil)$ be a $C^*$-algebra and $k\in \Proc(AtoAk$.
%Then the following are equivalent:\\
%{\rm (i)} $k$ is a Markov-regular, completely positive and contractive
%standard QS cocycle\\ \hspace*{0.5cm} on  $\alg$; \\
%{\rm (ii)}
%$k=k^{\phi}$ where $\phi\in CB\big(\alg ; M_{\khat}(\alg)\big)$
%satisfies $\phi(1) \leq 0$ and may be decomposed as follows:
%\begin{equation}
%\phi(a) =
%\Psi (a) -  a \ot \QSproj - \big(a\ot |e_0\ra\big) J - J^* \big(a \ot \la e_0|\big)
%\;\;\; a\in \alg, \label{CPC0}\end{equation}
%for a completely positive map $\Psi : \alg\to\alg''\,\overline{\ot} B(\kilhat)$ and
%operator $J\in \alg''\,\overline{\ot}\la\kilhat |$.
%\end{thm}
%Hee $\alg''$ denotes the double commutant of $\alg$ in $B(\hil)$.
%

\section{ $C^*$-hyperbialgebras and QS convolution cocycles}

In this section we describe the standard construction of new $C^*$-hyperbialgebras
via noncommutative conditional expectation, and give the definition
and some properties of quantum stochastic convolution cocycles. We then relate
quantum stochastic convolution cocycles on the respective $C^*$-hyperbialgebras.

\subsection*{$C^*$-hyperbialgebras and the conditional expectation construction}
\begin{defn}
A unital $C^*$-algebra $\alg$ is called a $C^*$-hyperbialgebra if it is
equipped with a unital completely positive map $\Com: \alg \to \alg
\otimes \alg$ (called a coproduct) and a character $\Cou:\alg \to \bc$
(called a counit) satisfying the following conditions:
\[ (\Com \ot \text{id}_{\alg}) \circ \Com =
(\text{id}_{\alg} \ot \Com) \circ \Com, \]
\[ (\Cou \ot \text{id}_{\alg}) \circ \Com =
(\text{id}_{\alg} \ot \Cou) \circ \Com  =
\text{id}_{\alg}.\]
If additionally $\Com$ is multiplicative then $\alg$ is called a
$C^*$-bialgebra (and may be thought of as a quantum compact semigroup with a neutral element).
\end{defn}

The following construction, of new $C^*$-hyperbialgebras from old, was
described explicitly (in the context of compact quantum hypergroups) in the
papers~\cite{Kal1} and~\cite{Kal2}, but its origins go back much further
(see~\cite{ChaV} and references therein). All known examples of
noncommutative $C^*$-hyperbialgebras arise in this way from
$C^*$-bialgebras.

\begin{propn}   \label{expecthyper}
Let $(\alg, \Com, \Cou)$ be a $C^*$-hyperbialgebra. Assume that
$\wt{\alg}$ is a unital $C^*$-subalgebra of $\alg$ and that there exists a
conditional expectation, that is a norm-one
projection, $P$ from $\alg$  onto $\wt{\alg}$ satisfying the following
identities:
\[
(P \ot \id_{\alg})\circ \Com\circ P =
(P \ot P)\circ \Com =
(\id_{\alg} \ot P)\circ \Com.
\]
Then $(\wt{\alg}, \wt{\Com}, \wt{\Cou})$ is a $C^*$-hyperbialgebra,
where
\[
\wt{\Com} = (P \ot P)\circ \Com|_{\wt{\alg}} \text{ and }
 \wt{\Cou} = \Cou|_{\wt{\alg}}.\]
\end{propn}

Two particular cases of this construction are \emph{double coset bialgebras}
and \emph{Delsarte $C^*$-hyperbialgebras}; they are described below.

Let $(\alg_1,\Com_1,\Cou_1)$ and $(\alg_2,\Com_2,\Cou_2)$ be
$C^*$-bialgebras and assume that the latter is a \emph{quantum
subsemigroup} of the former. This means that there exists a unital
*-homomorphism $\pi:\alg_1 \to \alg_2$ which is surjective and
intertwines the coalgebraic structures:
\[ (\pi \ot \pi)\circ \Com_1 = \Com_2 \circ \pi, \;\; \Cou_2 \circ \pi = \Cou_1.\]
Assume additionally that $\alg_2$ admits a Haar state; this means that
there exists a state $\mu \in \alg_2^*$ such that
for all $a \in \alg_2$
\[ (\mu \ot \id_{\alg_2}) \big(\Com_2 (a)\big) = \mu(a) 1_{\alg_2}.\]
 Define the following $C^*$-subalgebras of $\alg_1$:
\[\alg_1 /\alg_2 := \{ a \in \alg_1: (\id_{\alg_1 } \ot \pi ) \circ \Com_1
(a) = a \ot 1\},\]
\[ \alg_1 \backslash \alg_2 := \{ a \in \alg_1: (\pi \ot \id_{\alg_1 }  )
\circ \Com_1 (a) = 1 \ot a\}\]
\[ \alg_2 \backslash \alg_1 / \alg_2 :=
\alg_1 /\alg_2 \, \cap \, \alg_1 \backslash \alg_2.\]
called respectively the
\emph{algebras of left and right cosets of $\alg_2$}
and the \emph{double coset bialgebra}.

It can be checked that the map $P:\alg_1 \to \wt{\alg}:=\alg_2
\backslash \alg_1 / \alg_2$ defined by
\[ P(a) = \left((\mu \circ \pi) \ot \id_{\alg_1} \ot (\mu \circ \pi) \right) \left(\Com_1 \ot \id_{\alg_1} \right)
 \Com_1 (a), \;\;\;  a \in \alg_1,\]
satisfies the conditions given in Proposition~\ref{expecthyper}. Its
action may be understood as averaging (twice)
over the quantum subsemigroup; this construction is common in the theory of classical hypergroups (\cite{hyp}).

Let now $(\alg,\Com,\Cou)$ be a $C^*$-bialgebra and assume that a compact group $\Gamma$ acts
(continuously with respect to the topology of pointwise convergence) on $\alg$ by $C^*$-algebra automorphisms satisfying
\[ (\gamma \ot \gamma) \circ \Com = \Com \circ \gamma,\;\; \Cou \circ \gamma = \Cou, \;\;\; \gamma \in \Gamma.\]
Let $\wt{\alg}$ be the fixed point subalgebra,
$\wt{\alg}:=\{ a \in \alg: \forall_{\gamma \in \Gamma}\; \gamma(a) = a\}$.
It is easily checked that the map $P:\alg \to \wt{\alg}$ given by
\[ P(a) = \int_{\Gamma} \gamma(a) d\gamma, \;\;\; a \in \alg,\]
where $d\gamma$ denotes the normalised Haar measure on $\Gamma$, satisfies
the assumptions of Proposition~\ref{expecthyper}. The resulting
$C^*$-hyperbialgebra is called a \emph{Delsarte $C^*$-hyperbialgebra}.

Given a $C^*$-hyperbialgebra $\alg$, each operator space $\Vil$
determines a map
\begin{equation} \label{Rmap} R_{\Vil} : CB(\alg;\Vil) \to CB(\alg ;\alg \ot \Vil),
     \;\; \varphi \mapsto (\id_{\alg} \ot \varphi) \circ \Com. \end{equation}
%\[E_{\Vil} :CB(\alg;\alg \ot \Vil) \to CB(\alg ;\Vil),
%     \;\; \phi \mapsto (\Cou \ot \id_{\Vil} )\circ \phi .\]
When $\Vil = \bc$ we write $R$ instead of $R_{\bc}$.

\subsection*{QS convolution cocycles}

The following definition originates in \cite{LSqscc1} and is inspired by
the theory of quantum L\'evy processes. Let $\alg$ be a
$C^*$-hyperbialgebra.

\begin{defn}
A \emph{QS convolution cocycle}
on $\alg$, with noise dimension space $\kil$, is a completely bounded process
$l\in\ProcAtoC$ such that, for $s,t \geq 0$,
\begin{equation*}
%\label{inc}
l_{s+t} = (l_s \ot (\sigma_s \circ l_t)) \circ \Com
\text{ and }
l_0 = \iota_{\FFock}\circ \Cou.
\end{equation*}
The first of these conditions is referred to as the
\emph{convolution increment property}.
\end{defn}

A QS convolution cocycle $l$ is said to be \emph{Markov-regular}
if its \emph{Markov convolution semigroup of functionals}
$\lambda: \Rplus\to\alg^*$, defined by
\[ \lambda_t(a) =
\la \ve(0), l_t(a) \ve(0) \ra, \;\;\; t\geq 0, a \in \alg,\]
is norm-continuous.

For $\varphi\in CB\left(\alg;B(\kilhat)\right)$ we consider
\emph{coalgebraic
QS differential equations} of the form
\begin{equation}  \label{coalgQSDE}
dl_t =
l_t \star d \Lambda_{\varphi}(t) =
(l_t \star_{\pi} \varphi) d \Lambda_t
, \;\;\; l_0 = \iota_{\FFock}\circ \Cou,
\end{equation}
where $\pi$ indicates a tensor flip exchanging the order of $\kilhat$ and
$\Fock$. In fact the above equation may also be written as an equation of
the type (\ref{QSDE}), with $\phi= R_{B(\khat)} \varphi$ and
$\theta=\Cou$.
The unique solution of (\ref{coalgQSDE}) will be denoted by $l^{\varphi}$.
The process $l^{\varphi} \in \ProcAtoC$ is a Markov-regular
\emph{weak QS convolution cocycle}.
For full discussion of the precise meaning of (\ref{coalgQSDE}), weak QS
convolution cocycles and relations between the equation (\ref{coalgQSDE})
and equations of the type (\ref{QSDE}) we refer to~\cite{LSqscc2} and
\cite{LSqsde}.
The next two propositions are proved in \cite{LSqscc2} by applying the
semigroup decompositions of cocycles and convolution cocycles.

\begin{propn} \label{gencor}
Let $l = l^\varphi$ and $k = k^\phi$ where
$\varphi \in CB\big( \alg ; B(\kilhat) \big)$
and $\phi = R_{B(\kilhat)} \varphi$. Then the process $l$ is completely
bounded (respectively, completely positive and contractive) if and only if
$k$ is, and if $l$ and $k$ are completely bounded then
\begin{equation} \label{k-Rl}
k_t = R_{B(\FFock)} l_t, \;\; t \in \br_+ .
\end{equation}
\end{propn}

\begin{propn} \label{cocyclecor}
Let $k=R_{B(\FFock)} l$ where
$l$ is a completely bounded process in $\ProcAtoC$. Then
$l$ is a QS convolution cocycle if and only if
$k$ is a standard QS cocycle on $\alg$,
and in this case $l$ is Markov-regular if and only if $k$ is.
\end{propn}

Application of these results to the characterisation of the generators of
Markov-regular CPC QS cocycles, and *-homomorphic QS cocycles, given
in~\cite{lwblms} and~\cite{lwjlms} respectively, leads to the following
results.

\begin{thm} [\cite{LSqscc2}]\label{genstruct}
Let $\alg$ be a $C^*$-hyperbialgebra and $l\in \PProcAtoCk$.
Then the following are equivalent:\\
{\rm (i)} $l$ is a Markov-regular, completely positive and contractive QS
convolution cocycle;\\
{\rm (ii)}
$l=l^{\varphi}$ where $\varphi\in CB\big( \alg ; B(\khat)\big)$
satisfies $\varphi(1) \leq 0$ and may be decomposed as follows:
\begin{equation} \label{str1}
\varphi(a) =
\psi (a) - \Cou(a) \left(\QSproj+ |e_0\ra\la\chi | + |\chi\ra\la e_0| \right),
\;\;\; a\in \alg,
\end{equation}
for some completely positive map $\psi: \alg \to B(\khat)$ and vector $\chi\in\khat$.
\end{thm}

\begin{thm}[\cite{LSqscc2}] \label{mult}
Let $\alg$ be a $C^*$-bialgebra and let $l=l^\varphi$ where
$\varphi\in CB\big(\alg; B(\kilhat)\big)$. Then the following are
equivalent\textup{;} \\
{\rm (i)} $l$ is *-homomorphic\textup{;} \\
{\rm (ii)} $\varphi$ satisfies
\begin{equation} \label{homold}
\varphi(a^*b) =
\Cou(a)^* \varphi(b) +\Cou(b) \varphi(a)^* + \varphi(a)^* \QSproj \varphi(b),
\;\; \; a,b \in \alg.
\end{equation}
\end{thm}

% These are related in Proposition~\ref{hom} below.

\subsection*{QS convolution cocycles and the conditional expectation construction}
We end this section by describing the connection between QS
convolution cocycles on $C^*$-hyperbialgebras related by the construction given in
Proposition~\ref{expecthyper}.

\begin{propn} \label{corresp}
Let $(\wt{\alg}, \wt{\Com}, \wt{\Cou})$ be the $C^*$-hyperbialgebra
arising from a $C^*$-hyperbialgebra  $(\alg, \Com, \Cou)$ via the
construction presented in Proposition~\ref{expecthyper}, with associated
conditional expectation $P$.
Then there is a 1-1 correspondence between QS convolution cocycles on
$\wt{\alg}$ and $P$-invariant  processes on $\alg$ satisfying the
convolution increment property and having initial condition given by the
functional $\Cou \circ P$.
\end{propn}

\begin{proof}
Assume first that  $\wt{l}\in \AwtProcwtAtoC$ is a QS convolution cocycle
and define $l \in \ProcAtoC$ by
\[ l_t = \wt{l}_t \circ P, \;\; t \geq 0.\]
Then clearly $l_0(a) = \Cou (P(a))$ for all $a \in \alg$, and $l$ is $P$-invariant. It remains to
check it is a convolution increment process. Choose $t,s \geq 0$ and compute:
\begin{align*}
l_{s+t} =
\wt{l}_{s+t} \circ P &=
\left( \wt{l}_s\ot (\sigma_s\circ\wt{l}_t) \right)\circ\wt{\Com}\circ P \\
&=
\left( \wt{l}_s \ot (\sigma_s \circ \wt{l}_t) \right)
\circ(P \ot P)\circ \Com P \\
&=\left(\wt{l}_s\ot (\sigma_s\circ\wt{l}_t)\right)\circ (P\ot P)\circ\Com
=  \left( l_s \ot (\sigma_s \circ l_t) \right)\circ \Com.
\end{align*}

Conversely, let $l \in \ProcAtoC$ be a $P$-invariant convolution
increment process, with initial condition given by $\Cou \circ P$. Then
the process $\wt{l}\in \AwtProcwtAtoC$, defined simply by
the restriction of $l$, is a QS convolution cocycle on $\wt{\alg}$ ---
again the only thing to be checked
is the convolution increment property: for all $s,t \geq 0$, $a \in \wt{\alg}$,
\begin{align*}
\wt{l}_{s+t} (a) &= l_{s+t} (a) \\
&= \left( l_s \ot (\sigma_s \circ l_t) \right) \big(\Com (a)\big) \\
&=
\left( l_s \ot (\sigma_s \circ l_t) \right) \big(\Com (Pa)\big) \\
&=
\left(l_s\ot (\sigma_s\circ l_t)\right)
\big((P\ot P)\circ\Com\circ P\big)(a)
=
\left(\wt{l}_s\ot (\sigma_s\circ\wt{l}_t)\right)\big(\wt{\Com}(a)\big).
\end{align*}
\end{proof}
\noindent
Markov-regularity is clearly preserved in the above correspondence.

If $\Cou = \Cou \circ P$ (as is the case for Delsarte
$C^*$-hyperbialgebras, but usually not for double coset bialgebras), then
the processes $l$ arising in the proof of the above theorem are obviously
QS convolution cocycles. Assuming this is the case, it is easily checked
that if $\varphi\in CB(\wt{\alg}; B(\kilhat))$ then
$\wt{l}:=l^{\varphi} \in \AwtProcwtAtoC$ corresponds to the
process $l^{\varphi\circ P} \in \ProcAtoC$.
There is an analogous correspondence on the level of weak QS convolution
cocycles.

In \cite{LSqscc2} a variety of examples of $C^*$-bialgebras is presented,
and $^*$-homomorphic QS convolution cocycles on them are given alternative
interpretations.

\section{Generator of CPC QS convolution cocycle}

In this section we consider the detailed structure of the stochastic
generators of CPC QS convolution cocycles. Our approach is direct,
following ideas used in the study of CPC standard QS cocycles.
The crucial analysis was carried out in~\cite{lp}, with extension to
infinite dimensional noise done in \cite{lw1} and \cite{lwblms} (see also
\cite{Slava}).

Adapting arguments used in \cite{lp} requires some care, and the $R$-map
introduced in (\ref{Rmap}) is an indispensable tool. A
straightforward approach to complete positivity for QS convolution
cocycles leads to nontrivial considerations of the proper
convolution-counterpart of conditional complete positivity, and here the
$R$-map does not appear to be helpful.
However, nonnegative-definite kernels taking values in a
$C^*$-algebra
do behave well under the $R$-map, as will be seen in the proof of the
next proposition.

For the rest of this section $\alg$ denotes a fixed $C^*$-hyperbialgebra.
For any $\tau\in B(\alg)$ define $\partial{\tau}:\alg \times \alg \to \alg$ by
\[ \partial \tau (a_1,a_2) = \tau (a_1^*a_2) - a_1^* \tau(a_2) -
\tau(a_1^*) a_2 +
a_1^* \tau(1) a_2, \;\; a_1,a_2 \in \alg.\]
By analogy,  for any $f \in \alg^*$ define $\partial_{\Cou}f: \alg \times
\alg \to \bc$
by
\[ \partial_{\Cou}f (a_1, a_2) = f(a_1^* a_2) - \Cou(a_1^*) f(a_2) - f(a_1^*) \Cou(a_2) +
\Cou(a_1^*) f(1) \Cou(a_2), \;\; a_1,a_2 \in \alg.\]
If $(\rho,\Kil)$ is a representation of $\alg$, a map
$\delta: \alg \to B(\bc;\Kil)$
is called a \emph{($\rho,\Cou$)-derivation} if for all $a_1,a_2 \in \alg$
\[ \delta(a_1 a_2) = \rho(a_1) \delta(a_2) + \delta(a_1) \Cou(a_2).\]
Observe that if $\varphi \in CB(\alg;B(\kilhat)$ and $l^{\varphi}$ is completely positive then it is easily verified that it
is contractive too if and only if
\begin{equation}
\label{cont} \varphi(1)\leq 0.
\end{equation}

We need to start with the finite-dimensional situation.
The key fact is the following result, corresponding to Theorem 4.1 in
\cite{lp}.

\begin{lemma}  \label{findim}
Assume that $\kil$ is finite dimensional. Let $\varphi\in CB\big(\alg; B(\kilhat)\big)$ and suppose that the
\textup{(}weak\textup{)} QS convolution cocycle $l:=l^{\varphi}\in
\CdProcAtoCd$ is CPC. Then
there exist a unital *-representation $(\rho, \Kil)$ of $\alg$, a
($\rho,\Cou$)-derivation $\delta:\alg\to B(\bc;\Kil)$,
an operator $D \in B(\kil; \Kil)$ and a vector $d\in\kil$ such that
\begin{equation} \label{macierz00}
\varphi(a)=\begin{bmatrix}\lambda(a)&
\Cou(a)\langle d|+\delta^\dagger(a)D\\ \Cou(a)|d\ra +D^*\delta(a)&
D^*\rho(a)D-\Cou(a)I_{\kil}\end{bmatrix},
\;\;\; a\in \alg,
\end{equation}%**
where the functional $\lambda$ is real,
\[ \partial_{\Cou}\lambda (a_1, a_2) = \delta(a_1)^*\delta(a_2),
\;\; a_1,a_2 \in \alg,\]
and the following minimality condition holds:
\begin{equation} \label{minimality}
\Kil= {\clLin} \{\delta(a)1 + \rho(a) D c : a \in \alg, c \in \kil\}.
\end{equation}
If $(\Kil', \rho', \delta', D')$ is another quadruple satisfying the
above conditions \textup{(}except possibly the minimality
condition\textup{)}, then there exists a unique isometry $V: \Kil \to
\Kil'$ such that
\[ \delta'(a) = V \delta(a), \;\; \rho'(a) V = V \rho(a), \;\;D' = VD, \;\; a \in \alg.\]
\end{lemma}

\begin{proof}   The proof is a modification of the argument used in the
proof of Lemma 4.5 in \cite{lp}, where $\kil$ is taken to be $\bc^d$.
Write $\varphi$ in block matrix form:
\[
\left[\begin{matrix} \lambda
 & \widetilde{\eta} \\ \eta
 & \sigma - \Cou(\cdot)I_{\kil} \end{matrix} \right].
\]
By Propositions~\ref{gencor} and~\ref{cocyclecor},
$k=R_{B(\FFock_{\kil})} l$ is
a CPC standard QS cocycle and  $\phi=R_{B(\kilhat)} {\varphi}$ is real.
Therefore $\varphi$ is real too, in particular $\widetilde{\eta} = \eta^{\dag}$,
and $\phi$ has block matrix form
\[
\left[ \begin{matrix}{\tau}  & \alpha^{\dag} \\ \alpha
 & \nu - \iota \end{matrix} \right],
\]
where $\tau=R\lambda$, $\alpha=R_{B(\bc;\kil)}{\eta}$ and
$\nu=R_{B(\kil)}{\sigma}$. Now Lemma 4.4 in \cite{lp}
implies that the map $\Psi$ from $\alg\times\alg$ to $\alg\ot B(\kilhat)$,
there identified with $M_{d+1}(\alg)$, defined by
\[ \Psi(a_1, a_2) = \left[ \begin{matrix} \partial\tau (a_1,a_2)  & \alpha^{\dag} (a_1^*a_2) -
a_1^* \alpha^{\dag}(a_2)\\ \alpha (a_1^* a_2) - \alpha(a_1^*)a_2 & \nu(a_1^* a_2) \end{matrix}
\right], \;\; a_1,a_2 \in \alg, \]
is nonnegative-definite. Observe that if $\psi: \alg \times \alg \to
B(\kilhat)$ is defined by the formula
 \[ \psi(a_1, a_2) = \left[ \begin{matrix} \partial_{\Cou}\lambda (a_1,a_2)  & \eta^{\dag} (a_1^*a_2) -
\Cou(a_1^*) \eta^{\dag}(a_2)\\ \eta (a_1^* a_2) - \eta(a_1^*)\Cou(a_2) & \sigma(a_1^* a_2) \end{matrix}
\right], \;\; a_1,a_2 \in \alg, \]
then $ \psi = (\Cou \ot {\rm id}_{B(\kilhat)} ) \circ \Psi$.
This in turn implies that $\psi$
is a nonnegative-definite kernel. Indeed, for any
$n\in\bn,\, a_1, \ldots ,a_n\in\alg$ and $T_1, \ldots ,T_n\in B(\kilhat)$
\[ \sum_{i,j=1}^n T_i^* \psi(a_i, a_j) T_j =
\big(\Cou \ot {\rm id}_{B(\kilhat)}\big)
\Big(\sum_{i,j=1}^n (1_{\alg} \ot T_i^*) \Psi(a_i, a_j) (1_{\alg} \ot T_j) \Big) \geq 0, \]
as $(1_{\alg} \ot T_i)^* =  (1_{\alg} \ot T_i^*) \in \alg \ot B(\kilhat)$,
$\Cou$ is CP, and $\Psi$ is nonnegative-definite.

Now let $(\Kil, \chi)$ be the minimal Kolmogorov pair associated with $\psi$. This
means that $\chi$ is a map $\alg \to B(\kilhat; \Kil)$ satisfying
\begin{align*}
& \chi(a_1)^* \chi(a_2) = \psi(a_1, a_2), \;\;\;\; a_1, a_2 \in \alg, \\
& \Kil = {\clLin}\{\chi(a) \zeta: a\in \alg, \zeta \in \kilhat\}.
\end{align*}
Properties of $\psi$ imply that $\chi$ is linear and bounded. Write
$\chi = [\delta \; \gamma]$, where $\delta\in B\big(\alg;B(\bc;\Kil)\big)$
and $\gamma \in B(\alg; B(\kil;\Kil))$. Then, for any $a,b\in \alg$,
\[ \delta(a)^* \delta(b) = \partial_{\Cou}\lambda (a,b) \text{ and }
\gamma (a)^*\delta(b) = \eta(a^*b) - \eta(a^*)\Cou(b). \]
Setting $a=b=1$ shows that $\delta(1)=0$.
Now for $u \in \alg$ unitary, define
\[
\delta_u (a) = \delta(ua) - \delta(u) \Cou(a), \
\gamma_u(a) = \gamma(ua) \text{ and } \chi_u = [\delta_u \; \gamma_u],
\quad \text{ for } a\in\alg.
\]
A straightforward computation yields
\[ \chi_u (a_1)^* \chi_u (a_2) =   \chi (a_1)^* \chi (a_2).\]
The uniqueness of minimal Kolmogorov pairs implies the existence of a unique
isometry $\rho(u): \Kil \to \Kil$
given by the formula
\[ \rho(u) (\delta(a)1 + \gamma(a)c) =
\delta(ua)1 - \delta(u) \Cou(a) + \gamma(ua)c,
\;\; a \in \alg, c \in \kil.\]
It follows, by standard arguments, that
\[ \rho(a) (\delta(b)1 + \gamma(b)c) =  \delta(ab)1 - \delta(a) \Cou(b) + \gamma(ab)c,
\;\; a,b \in \alg, c \in \kil,\]
defines a bounded operator $\rho(a)$ on $\Kil$. Moreover, it is easily checked that the resulting map
$\rho: \alg \to B(\Kil)$
is indeed a *-representation of $\alg$. It immediately follows that $\delta$ is a
($\rho,\Cou$)-derivation and also, from minimality and the identity $\delta (1)=0$, that
$\rho$ is unital.
Put $D= \gamma(1)\in B(\kil;\Kil)$. Then $\gamma(a) = \rho(a)D$, and
furthermore
$\sigma(a) = D^* \rho(a)D$ and
$\eta(a) = \Cou(a)\eta(1) + D^* \delta(a)1$. This yields
($\ref{macierz00}$) with $d=\eta (1)1$.

The second part of the lemma follows once more from uniqueness of the
Kolmogorov construction.
\end{proof}

The step from finite-dimensional to arbitrary noise dimension space follows
in exactly the same way as for standard cocycles.

\begin{lemma}  \label{infdim}
Assume that $\kil$ is an arbitrary Hilbert space. Let $\varphi\in CB(\alg; B(\khat))$ and suppose that the
\textup{(}weak\textup{)} QS convolution cocycle $l^{\varphi}\in
\PProcAtoCk$ is CPC. Then the conclusions of Lemma~\ref{findim} hold.
%without the assumption that $\kil$ is finite dimensional.
%Then there exist a unital *-representation
%$(\rho,\Kil)$ of $\alg$, a ($\rho,\Cou$)-derivation
%$\delta:\alg\to B(\bc;\Kil)$, an operator $D \in B(\kil; \Kil)$ and a
%vector $d\in \kil$ such that
%\begin{equation*}
%\varphi(a)=\begin{bmatrix}\lambda(a)&\Cou(a)\langle
%d|+\delta^\dagger(a)D\\
%\Cou(a)|d\ra +D^*\delta(a)&D^*\rho(a)D-\Cou(a)I_{\kil}\end{bmatrix},
%\end{equation*}
%where $\lambda$ is real,
%\[ \partial_{\Cou}\lambda (a_1, a_2) =  \delta(a_1)^* \delta(a_2),
%\;\; a_1,a_2 \in \alg,\]
%and the minimality condition~\eqref{minimality} holds.
\end{lemma}

\begin{proof}
Let $\{\kil_\iota:\iota \in \mathcal{I} \}$ be an indexing of the set of all
finite-dimensional subspaces of $\kil$, which is partially ordered by inclusion.
As in \cite{lw1} we consider finite-dimensional cut-offs of  both $l^{\varphi}$ and
$\varphi$ itself. For each $\iota \in\mathcal{I}$ denote by $\varphi_\iota$
the map $\alg \to B(\widehat{\kil_\iota})$ given by the formula
\[ \varphi_\iota(a) = P_\iota \varphi(a)  P_\iota, \;\; a \in \alg,\]
where $P_\iota\in B(\kilhat)$ is the orthogonal projection onto $\kilhat_\iota$.
Setting $l^{(\iota)} = l^{\varphi_\iota}$, $\FFock^\iota = \FFock^{\kil_\iota}$,
$\Exps^{\iota}=\Exps_{\kil_{\iota}}$ and letting $\mathbb{E}_\iota$ denote the vacuum conditional expectation
map from $B(\FFock_{\kil})$ to $B(\FFock^\iota)$,
it is easy to see that $l^{(\iota)} \in \iotProcAtoCiota$ is a
CPC QS convolution cocycle and that it satisfies
\[
l^{(\iota)}_t(a) = \mathbb{E}_\iota [l^{\varphi}_t(a)], \;\;
a\in \alg, t\in\Rplus.
\]

Lemma \ref{findim} yields quadruples
$(\Kil_\iota, \rho_\iota, \delta_\iota, D_\iota)$, unique up to isometric
isomorphism, such that for all $a \in \alg$
\begin{equation*}
\varphi_\iota(a) =
\begin{bmatrix}\lambda(a)&\Cou(a)\langle d_\iota|+\delta_\iota^\dagger(a)D_\iota\\
\Cou(a)|d_\iota\ra +D_\iota^*\delta_\iota(a)&D_\iota^*\rho_\iota(a)D_\iota-\Cou(a)I_\iota\end{bmatrix},
%\label{macierz3}
 \end{equation*}
where $I_\iota$ denotes the identity operator on $\kil_\iota$.

Exploiting uniqueness one can construct an inductive limit $\Kil$
of the Hilbert spaces $\Kil_\iota$. Denote by $U_\iota$ the respective isometry
$\Kil_\iota\to\Kil$. Then there is a unital *-representation $\rho$
of $\alg$ on $\Kil$, a ($\rho,\Cou$)-derivation
$\delta: \alg \to B(\bc;\Kil)$ and, for each $c \in \kil$ a vector $c_D \in \Kil$ such that
\[\rho(a)U_\iota = {U_\iota}\rho_\iota(a),
\ \
\delta(a) = U_\iota\delta_\iota(a)
\text{ and }
{c_D}  = U_\iota D_\iota c,\]
for all $\iota \in \mathcal{I}$, $a \in \alg$ and $c \in \kil_\iota$.
The map $c \mapsto c_D$ is linear; it remains to show that it is bounded.
To this end
observe that, for any $\iota\in\mathcal{I}$ such that $c\in\kil_\iota$,
\begin{align*}
\Big\la \binom{0}{c}, \varphi(1) \binom{0}{c} \Big\ra &=
\Big\la \binom{0}{c}, \varphi_\iota(1) \binom{0}{c} \Big\ra \\ &=
\big\la c, (D_\iota{^*} D_\iota - \Cou(1) I_\iota )c \big\ra =
\|D_\iota c \|^2 - \|c\|^2 = \|c_D\|^2 - \| c\|^2,
\end{align*}
and inequality (\ref{cont}), being a consequence of the contractivity
of $l^{\varphi}$, implies that $\|c_D\| \leq \|c\|$. The operator
$D\in B(\kil; \Kil)$ given by $Dc = c_D$ completes the tuple whose existence
we wished to establish. Minimality holds by construction.
\end{proof}

Automatic innerness of ($\rho,\Cou$)-derivations (see~\cite{LSqscc2})
leads to the following theorem.

\begin{thm}  \label{spec}
Let $\varphi\in CB(\alg;B(\khat))$, for a $C$*-hyperbialgebra $\alg$, and
suppose that the weak QS convolution cocycle
$l^{\varphi}\in \PProcAtoCk$ is completely positive and contractive. Then
there exists a tuple $(\Kil,\rho,D,\xi,d,e,t)$ constisting of a
unital *-representation $(\rho,\Kil)$ of $\alg$, a contraction $D\in
B(\kil;\Kil)$, vectors $\xi\in\Kil$ and $d,e \in \kil$, and a real number
$t$, such that
\begin{equation}  \label{specified}
\varphi(a)=\begin{bmatrix}\lambda(a)&\Cou(a)\langle d|+\delta^\dagger(a)D\\
\Cou(a)|d\ra +D^*\delta(a)&D^*\rho(a)D-\Cou(a)I_{\kil}\end{bmatrix},
\end{equation}
$t=\lambda(1)\leq 0$ and $d=(I_{\kil}-D^*D)^{1/2}e$, $\|e\|^2\leq -t$, and,
for all $a \in \alg$,
\begin{equation} \label{deltalambda}
\delta(a) = \big(\rho(a) - \Cou(a)\big)|\xi\ra, \;\;\;
\lambda(a) =
\Cou(a) (t - \|\xi\|^2)+ \langle \xi, \rho(a) \xi \rangle.
\end{equation}
\end{thm}

\begin{proof}
Lemma \ref{infdim} gives the form~\eqref{specified} for some
$\rho, \Kil, \delta$ and $D$.
As all ($\rho,\Cou$)-derivations are inner (see \cite{LSqscc2}, Appendix),
there exists $\xi \in \Kil$ such that
\[ \delta(a) = \rho(a) |\xi\ra - \Cou(a) |\xi\ra.\]

It remains to note that
\begin{equation}  \varphi(1)=
\left[ \begin{matrix} t &  \langle d| \\
|d \rangle & D^*D - I_{\kil}\end{matrix}\right],\label{contrac}\end{equation}
and the condition $\varphi(1)\leq 0$ implies contractivity of $D$, negativity of $t$ and the existence
of a vector $e\in \kil$ satisfying all the conditions above
(see the characterisation of positive matrices given in Lemma 2.1 of \cite{dilate}).
\end{proof}

\begin{remarks}
The converse is also true --- if a map
$\varphi\in CB\big(\alg;B(\kil)\big)$ has the above form then $l^\varphi$
is CPC. This follows easily from Propositions \ref{gencor} and \ref{cocyclecor} and Theorem 2.4 of
\cite{lwblms}.

If $l^{\varphi}$ is also unital, then $\varphi(1)=0$
and~\eqref{specified} takes the form
\[
\varphi(a)=\begin{bmatrix}\lambda(a) & \delta^\dagger(a)D\\
D^*\delta(a) & D^*\rho(a)D-\Cou(a)I_{\kil}\end{bmatrix},
\]
with $D$ being an isometry. This corresponds exactly to the characterisation obtained
in the purely algebraic case by Franz and Sch\"urmann (\cite{UweSch}).

The characterisation in Theorem \ref{spec} yields, as announced in the beginning of this section,
an alternative proof of the existence of the decomposition \eqref{str1} of Theorem \ref{genstruct}.
Indeed, for $\varphi:\alg \to B(\kilhat)$ of the form (\ref{specified}),
define $S: \khat\to \Kil$ by $S = [ |\xi\ra \;\; D ]$. Then
\[ \varphi(a) = S^* \rho(a) S +
\left[ \begin{matrix}{\lambda_0(a)} & \Cou(a)\langle d - D^* c | \\
\Cou(a) |d - D^* c\rangle & - \Cou(a)I_{\kil} \end{matrix}\right], \;\; a \in \alg, \]
where $\lambda_0 (a) = \lambda(a) -\langle \xi, \rho(a) \xi \rangle$.
Note that as $\partial_{\Cou}\lambda_0 (a_1, a_2) = 0$
for any $a_1, a_2 \in \alg$, $\lambda_0= \lambda_0(1) \Cou$ -
one can check that $\big(\lambda_0 -  \lambda_0(1) \Cou\,\big)$ is an
$(\Cou,\Cou)$-derivation and so is zero.
The map $\psi: \alg \to B(\kilhat)$ defined by
\[ \psi(a) = S^* \rho(a) S,\;\;\;   a \in \alg,\]
is  completely positive.
Setting $\chi = \binom{\frac{1}{2}\lambda_0(1)}{d-D^*c}$
yields
the required decomposition.
\end{remarks}

\section{Dilations to *-homomorphic QS convolution cocycles}

This section addresses the question of dilating a completely positive,
contractive QS convolution cocycle on a $C^*$-bialgebra to a
*-homomorphic QS convolution cocycle. It is closely patterned on \cite{dilate}.
From now on we assume that $\alg$ is a $C^*$-bialgebra. Recall that this
means that $\alg$ is a $C^*$-hyperbialgebra whose  coproduct is
multiplicative. Let $\kil_0$ be a closed subspace of a noise dimension
space $\kil$. The vacuum conditional expectation from $B(\Fock_{\kil})$ to
$B(\Fock_{\kil_0})$ will be denoted by $\be_0$.

\begin{defn}
A (weak) QS convolution cocycle $j \in \PProcAtoCk$ is said to be a
stochastic dilation of a (weak) QS convolution cocycle $l\in
%\pPProcAtoCk0
\salvation$ if
\[ l_t = \be_0 \circ j_t, \;\;\; t \geq0.\]
\end{defn}

The following result follows in exactly the same way as its counterpart
for standard cocycles (\cite{dilate}, Lemma 1.2).

\begin{propn} \label{dilgen}
Let $\varphi\in CB(\alg;B(\khat))$ and $\psi\in CB(\alg; B(\khat_0))$, and
let  $j=l^{\varphi} \in \PProcAtoCk$ and $l=l^{\psi}\in
%\pPProcAtoCk0
\salvation$
be the respective QS convolution cocycles. Then $j$ is
a stochastic dilation of $l$ if and only if
$\psi (\cdot) = P_0 \varphi(\cdot) P_0$, where $P_0\in B(\khat)$ denotes
the orthogonal projection onto $\khat_0$.
\end{propn}

Generators of *-homomorphic cocycles may be described in the following
way.

\begin{propn} \label{hom}
Let $(\Kil, \rho, D, \xi, d, t)$ be a tuple as in Theorem \ref{spec}
and let $\varphi$ be the map in $CB\big(\alg;B(\khat)\big)$ given by the
formulas \eqref{specified} and \eqref{deltalambda}.
Then the \textup{(}weak\textup{)} QS convolution cocycle
$l^{\varphi}\in \PProcAtoCk$ is *-homomorphic if and only if the
following conditions hold\textup{:} \\
{\rm (i)} $D$ is a partial isometry,\\
{\rm (ii)} $Dd = 0$,\\
{\rm (iii)} $DD^* \in \rho(\alg)'$,\\
{\rm (iv)}  $t = - \| d \|^2,$\\
{\rm(v)} $DD^* \delta = \delta$.\\
where $\delta$ is the ($\rho,\Cou$)-derivation
$a\mapsto\big(\rho(a)-\Cou(a)I_{\kil}\big)|\xi\ra$.
\end{propn}

\begin{proof}
%By Theorem 6.1 and Proposition 4.4 of~\cite{LSqscc2},
%$j$ is *-homomorphic if and only if $\varphi$ is real and satisfies
%\begin{equation} \varphi(ab) = \varphi(a) \QSproj \varphi(b) + \Cou(a) \varphi(b) +
%\Cou(b) \varphi(a),\;\; \; a,b \in \alg.\label{homold}\end{equation}
In the language of Theorem $\ref{spec}$, the structure
relations~\eqref{homold} translate into the following identities:
\begin{align*}
&D^* \rho(a) DD^* \rho(b) D = D^* \rho(ab) D,  \\
&D^* \delta(ab) + \Cou(ab)|d\ra  =
D^* \rho(a) D \big( D^* \delta(b) + \Cou(b)|d\ra\big) + D^* \delta(a) \Cou(b) + \Cou(a)\Cou(b)|d\ra,  \\
&\lambda(a^*b) = \big\langle D^* \delta(a)1 + \Cou(a) d, D^* \delta(b)1 + \Cou(b) d \big\rangle +
\lambda(a^*)\Cou(b)
+ \Cou(a^*) \lambda(b),
\end{align*}
for all $a,b \in \alg$. As in Proposition 3.3 of \cite{dilate},
this in turn may be shown to be equivalent to the conditions (i)-(v).
\end{proof}

\begin{remarks}
Observe that the above characterisation excludes the possibility of obtaining exchange free dilations --- it can be
seen directly from ($\ref{homold}$) that if a Markov-regular *-homomorphic QS convolution cocycle is
generated by a map having the form
\[
\begin{bmatrix}*&*\\ *&0\end{bmatrix}
\]
then it is identically 0. This uses the fact that
$(\Cou,\Cou)$-derivations are trivial.
As to creation/annihilation free dilations, they are possible only for those CPC QS convolution cocycles, whose generators have the form
\[\begin{bmatrix}0&0\\ 0&*\end{bmatrix}.\]
%- the existence of the dilation follows then from the classical Stinespring theorem.

Moreover, $j$ is unital, as well as being *-homomorphic, if and only if
(iii) and (v) are satisfied, $D$ is an isometry, $d=0$ and $t=0$.
\end{remarks}

As a consequence of Theorem $\ref{spec}$ and Proposition~$\ref{hom}$, we
obtain the existence of stochastic dilations.

\begin{thm}
Every Markov-regular CPC QS convolution cocycle on a $C^*$-bialgebra $\alg$
%with separable noise dimension space,
admits a Markov-regular *-homomorphic stochastic dilation.
\end{thm}

\begin{proof}
%A simplified version of reasoning in proof of Theorem 4.1 of \cite{dilate} will do.
Let $l\in \AtoCProc$ be a Markov-regular CPC QS convolution cocycle.
Then $l=l^{\varphi}$
for some  $\varphi\in CB\big(\alg; B(\khat_0)\big)$ and we can assume that
$\varphi$ has matrix form
(\ref{specified}) for a tuple $(\Kil, \rho, D, \xi, d, e, t)$ with the
properties described
in Theorem \ref{spec}. Let $\kil_1, \kil_2$ be Hilbert spaces,
suppose that $d_1 \in \kil_1$, $d_2 \in \kil_2$, $D_1 \in B(\kil_1; \Kil)$
(all as yet unspecified) and consider the
map $\psi:\alg \to B(\kilhat)$, where $\kil :=
\kil_0 \oplus \kil_1 \oplus \kil_2$, given by ($a \in \alg$)

\begin{equation}  \psi(a) = \left[ \begin{matrix}{\lambda(a)}
 & \Cou(a)\langle d|+\delta^\dagger(a)D
 & \Cou(a)\langle d_1|+\delta^\dagger(a)D_1
 & \Cou(a)\langle d_2| \\
\Cou(a)|d\ra +D^*\delta(a) &
D^*\rho(a)D - \Cou(a)I_0
 & D^* \rho(a) D_1 & 0 \\
\Cou(a)|d_1\ra +D_1^*\delta(a) &
D_1^* \rho(a) D &
  D_1^*\rho(a)D_1 - \Cou(a)I_1 & 0 \\
\Cou(a)|d_2\ra  &
0 & 0& - \Cou(a)I_2
 \end{matrix}\right],\label{dilated} \end{equation}
with $I_i$ denoting $I_{\kil_i}$, $i=0,1,2$. Now
observe that $\psi$ can also be written in the form
\begin{equation}  \psi(a) = \left[ \begin{matrix}
{\lambda(a)} & \Cou(a)\langle \wt{d}|+\delta^\dagger(a)\wt{D}\\
\Cou(a)|\wt{d}\ra +\wt{D}^*\delta(a)&
\wt{D}^*\rho(a)\wt{D} - \Cou(a)I_{\kil} \end{matrix}\right],\label{specified2} \end{equation}
where
\[
\wt{d}=  \left( \begin{matrix} d \\ d_1 \\ d_2 \end{matrix}
\right) \in \kil \text{ and }
\wt{D} = \left[ \begin{matrix} D & D_1 & 0 \end{matrix} \right]
\in B(\kil;\Kil).
\]
As $\psi$ is clearly completely bounded, it generates a weak QS
convolution cocycle $l^{\psi}\in \PProcAtoCk$. It follows from
Proposition~\ref{dilgen} that $l^{\psi}$ is a stochastic dilation of
$l^{\varphi}$; it
remains to show that we can choose the parameters
$\kil_1, \kil_2$, $d_1, d_2$ and $D_1$ so that $l^{\psi}$ is
*-homomorphic.

To this end, it suffices to put $\kil_1 = \Kil$, $\kil_2 = \bc$,
\[
D_1 =  \big( I_1 - DD^*	  \big)^{\frac{1}{2}}, \; \;
d_1 = D e, \;\; d_2 = \sqrt{-(t + \|e\|^2)}.
\]
The above definitions make sense as $\|e \|\leq -t$ and $D$ is a
contraction.
%Put $ \wt{e} = \left( \begin{matrix} e \\ 0 \\ d_2 \end{matrix} \right)$.
%It is easily verified that the tuple
%$(\Kil, \rho, \xi, \wt{D}, \wt{d}, \wt{e}, t)$
%satisfies the conditions formulated in Theorem \ref{spec}.
It remains then
to check properties (i)-(v) of Proposition~\ref{hom}.
First note that
\[ \wt{D}\wt{D}^* = DD^* + I_1 - DD^* = I_1,\]
which implies that conditions (i), (iii) and (v)  are satisfied (one can easily check
that $\wt{D}^* \wt{D}$ is a selfadjoint projection). Further
we obtain (ii):
\[
\wt{D} \wt{d} =
D (I_0-D^*D)^{\frac{1}{2}} e +
\left( I_1 - DD^*\right)^{\frac{1}{2}}D e = 0.
\]
Finally (iv) follows since
\[
\|\wt{d}\|^2 =
\|(I_0-D^*D)^{1/2} e\|^2 + \| D e\|^2 - \big(t+\|e\|^2\big) = -t.
\]
This completes the proof.
\end{proof}

If $l$ is unital and $\dim \Kil = \dim \Ran(I_{\Kil} - DD^*)$, then it is possible
to obtain the unital *-homomorphic dilation $j \in \PProcAtoCk$ of $l$ (with the noise dimension space $\kil=\kil_0 \oplus \Kil$).

\section{Stinespring theorem for QS convolution cocycles}

As the previous section was a variation on the theme of \cite{dilate},
this one addresses the convolution counterpart of the problem considered in \cite{Stine} for
standard QS cocycles. We shall show (in Theorem~\ref{Stin}) that each
Markov-regular, completely positive, contractive QS convolution cocycle
has a Stinespring-like decomposition in terms of a *-homomorphic cocycle
perturbed by a contractive process.

First we need some remarks on QS differential equations of the type:
\begin{equation}
{\rm d}W_t = F_t \big(I_{\khat} \ot W_t\big)\, {\rm d}\Lambda_t, \;\;\; W_0 =
\idf,\label{eqgen}
\end{equation}
where $F\in \ProckilhatE$ is a bounded process.
We say that $W$ is a weak solution of the above equation if for all
$f,g \in \Step$ and $t \geq0$
\[ \la \ve(f), (W_t - \idf) \ve(g) \ra  =
\int_0^t {\la \fhat(s) \ot \ve(f),
F_s  (I_{\kilhat} \ot W_s) (\ghat(s) \ot \ve(g) \ra ds}.\]

The solution of the above equation is given by the iteration procedure:
\begin{align*}
& X_t^0 = \idf, \;
X_t^1 = \int_0^t F_s (I_{\kilhat} \ot X_s^0) d\Lambda_s,\;\cdots\;,
X_t^{n+1} = \int_0^t F_s (I_{\kilhat} \ot X_s^n)
d\Lambda_s,\;\cdots \\
& W_t \ve(f) := \sum_{n=0}^{\infty} X_t^n \ve (f).
\end{align*}
Sufficient conditions for the above heuristics to be justified are that
$F$ is strongly measurable and has locally uniform bounds;
this is also sufficient for the uniqueness of strongly regular strong
solutions of the equation (\cite{Stine}, Proposition 3.1).
These conditions are clearly satisfied when
\[F_s = (\id_{B(\kilhat)} \ot l_s)(T), \;\;\; s\geq 0,\]
where $l$ is a Markov-regular, CPC QS convolution cocycle and $T\in
B(\khat) \ot \alg$.

Now let $j$ be the *-homomorphic QS convolution cocycle $l^{\varphi}$
 ($\varphi\in CB (\alg; B(\khat))$) and let $T\in B(\khat)\ot \alg$.
Assume that $W\in \ProcE$
is a \emph{bounded} solution to the equation
\begin{equation} \label{geneq} {\rm d}W_t =
(\id_{B(\kilhat)} \ot
j_t) (T) \big(I_{\khat} \ot W_t\big) {\rm d} \Lambda_t,\;\;\; W_0 = \idf. \end{equation}
We shall identify sufficient conditions for $W$ to be a contractive
process later.
The next question to be addressed is: when can we expect a process $k\in
\ProcAtoC$ defined by
\[k_t(a) = j_t(a) W_t,\;\;\; a\in \alg, \; t\geq 0,\]
to be a Markov-regular QS convolution cocycle?

The quantum It\^o formula yields
\begin{align*}
\big\langle \ve(f), k_t (a) \ve(g) \big\rangle
= &
\big\langle j_t(a^*) \ve(f), W_t \ve (g) \big\rangle \\
= &
\Cou(a) \la \ve(f),\ve(g) \ra + \\
& \int_0^t \textrm{ds}
\Big( \big\la \tilde{j}_s (I_{\kilhat} \ot a^*)(\fhat(s) \ot
\ve(f)), \tilde{j}_s(T)   \wt{W}_s (\ghat(s) \ot \ve(g)) \big\ra + \\
& \quad \quad
\big\la \tilde{j}_s (\phi(a^*)) (\fhat(s) \ot \ve(f)),  \wt{W}_s (\ghat(s)
\ot \ve(g)) \big\ra  + \\
& \quad \quad
\big\la \tilde{j}_s (\phi(a^*)) (\fhat(s) \ot \ve(f)),
(\QSproj\ot \idf)\tilde{j}_s(T)\wt{W}_s (\ghat(s)\ot\ve(g))\big\ra\Big)
\end{align*}
($f,g \in \Step, t\geq0$), where $\phi = R_{B(\kilhat)}\varphi$,
$\tilde{j}_s = (\idB \ot j_s)$ and $\widetilde{W}_s = I_{\kilhat} \ot W_s$.
Defining analogously $\tilde{k}_s = (\idB \ot k_s)$ we see that the above
equation may be written as
\begin{align*}
&\big\langle \ve(f), k_t (a)  \ve(g) \big\rangle =
\Cou(a) \langle \ve(f), \ve(g) \rangle + \\
&
\int_0^t \textrm{ds} \left( \Big\la \fhat(s) \ot \ve(f),
\tilde{k}_s \left((I_{\kilhat} \ot a)T + \phi(a) +
\phi(a)(\QSproj\ot\ida)T\right) (\ghat(s) \ot \ve(g)) \Big\ra \right).
\end{align*}

 The process  $k$ is equal to $l^{\psi}$ for some  $\psi \in CB (\alg;
B(\khat))$ if and only if $\wt{\psi} := (\psi \ot \id_{\alg})\circ\Com$ is given by
\begin{equation} \label{eq1}
a\mapsto
(I_{\kilhat} \ot a)T + \phi(a) + \phi(a) (\QSproj\ot \ida)T.
\end{equation}
Note that we need to work with the \emph{left version} of the map $R$
introduced in \eqref{Rmap} because of the tensor flip
in the definition of the coalgebraic QS differential equation (\ref{coalgQSDE}).
Let $\tau = (\idB \ot \Cou ) ( T) \in B(\khat)$.
Then $(\ref{eq1})$ implies that
\begin{equation} \label{sigm} \psi(a) = \Cou(a) \tau + \varphi(a) (1 +
\QSproj\tau),\end{equation}
and so
\begin{equation} \label{eq2} \wt{\psi}(a) = \tau \ot a +  \phi(a) +
\phi(a) (\QSproj\tau \ot\ida).  \end{equation}
Comparing $(\ref{eq1})$ with $(\ref{eq2})$ yields
\begin{equation} (I_{\kilhat} \ot a)T + \phi(a) (\QSproj\ot \ida)T =  \tau
\ot a + \phi(a) (\QSproj\tau \ot\ida). \label{compar}\end{equation}
If $T = \tau \ot \ida$ then this condition is automatically satisfied. If
$j$ is unital, then $T= \tau \ot \ida$ is also necessary for \eqref{compar}
to hold: put $a=\ida$ and use $\phi(\ida)=0$.

Observe that when $T = \tau \ot \ida$ the equation $(\ref{geneq})$ takes
the simple form
\begin{equation} \label{pert} {\rm d}W_t = (\tau \otimes U_t W_t ){\rm d}
\Lambda_t,\;\;\; W_0 = \idf, \end{equation}
with $U_t = j_t(1)$.  In this case the condition on $\tau$ assuring
contractivity of $W$ is also particularly simple.

\begin{thm} \label{gener}
Let $j=l^{\varphi}$ where $\varphi\in CB (\alg; B(\khat))$ and $\alg$ is a
$C^*$-bialgebra. Suppose that $j$ is *-homomorphic and $\tau\in B(\khat)$ satisfies the condition
\begin{equation} \label{con} \tau + \tau^* +\tau^*\QSproj\tau \leq 0. \end{equation}
Then the equation (\ref{pert}), with $U_t:= j_t(1)$, has a
unique contractive strong solution $W\in \ProcE$
\textup{(}contractive means that each $W_t$ is a contraction\textup{)}.
Moreover the process $W_t^* j_t(\cdot) W_t$ is equal to $l^{\theta}$,
where
\[ \theta(a) = \Cou(a) \left( \tau^*+\tau + \tau^*\QSproj \tau \right) +
(I_{\kilhat} + \tau^* \QSproj)\varphi(a) (I_{\kilhat} + \QSproj\tau),
\;\;\; a\in \alg.\]
\end{thm}

\begin{proof}
The discussion before the theorem shows that the equation (\ref{pert}) has a unique
strongly regular strong solution $W\in \ProcE$. The It\^o formula
yields, for
$u = \sum_{i=1}^k {\lambda_i \ve(f_i)}$, $k\in \bn$,
$\lambda_1, \ldots, \lambda_k\in \bc$,
$f_1,\ldots, f_k\in \Step$,
\begin{align*}
\la W_t u, W_t u \ra - \la u, u \ra =
&
\sum_{i,j=1}^k \overline{\lambda_i} \lambda_j
\int_0^t \textrm{ds} \Big( \big\la \fhat_i(s) \ot \ve(f_i),
\tau\fhat_j(s) \ot U_s\ve(f_j) \big\ra \\
& \qquad \qquad +
\big\la \tau\fhat_i(s) \ot U_s\ve(f_i),  \fhat_j(s) \ot \ve(f_j) \big\ra  \\
& \qquad \qquad +
\big\la \tau\fhat_i(s) \ot U_s\ve(f_i),  \QSproj\tau\fhat_j(s) \ot
U_s\ve(f_j) \big\ra\Big).
\end{align*}
As $U_s=j_s(1)$ and $j$ is *-homomorphic, each $U_s$ is a projection.
Therefore putting
$x(s) = \sum_{i=1}^k {\lambda_i \fhat_i(s) \ot U_s \ve(f_i)}$, $s \in [0,t]$, yields
\[ \la W_t u, W_t u \ra - \la u, u \ra=  \int_0^t \textrm{ds} \left\la x(s), \left((\tau+\tau^*+
\tau^*\QSproj\tau) \ot \idf\right) x(s) \right\ra \leq 0.\]
It follows that $W$ is contractive.

The proof of the second part of the theorem is a combination of the
considerations before its formulation and one more application of the It\^o
formula. Again let $f,g\in \Step$, $t \geq 0$, $a \in \alg$ and
$T=\tau \ot 1_{\alg}$,
let $\tilde{j}$, $\tilde{k}$, $\wt{W}$ and $\psi$ be defined as in the
discussion before the theorem and set
$\wt{\psi} = (\psi\ot\id_{\alg})\circ\Com$. Then
\begin{align*}
\langle \ve(f), W_t^* j_t(a) W_t\ve (g) \rangle
= &
\langle W_t \ve(f), j_t(a) W_t \ve (g) \rangle \\
= &
 \Cou(a) \la \ve(f),\ve(g) \ra + \\
& \ \
\int_0^t \textrm{ds}
\Big( \big\la \widetilde{W}_s (\fhat(s) \ot \ve(f)),
\tilde{k}_s(\wt{\psi}(a))   (\ghat(s) \ot \ve(g)) \big\ra + \\
& \qquad
\big\la \tilde{j}_s (T) \widetilde{W}_s  (\fhat(s) \ot \ve(f)),
\tilde{j}_s(I_{\kilhat} \ot a)
   \widetilde{W}_s (\ghat(s) \ot \ve(g)) \big\ra + \\
& \qquad \big\la \tilde{j}_s (T) \widetilde{W}_s (\fhat(s) \ot \ve(f)),
(\QSproj\ot \idf) \tilde{k}_s(\wt{\psi}(a))(\ghat(s) \ot \ve(g))
\big\ra \Big).
\end{align*}
Finally,  (\ref{sigm}) yields
\begin{align*}
\big\la \ve(f), W_t^* j_t(a) W_t\ve(g) \big\ra
=&
\Cou(a) \big\la \ve(f),\ve(g) \big\ra \\
& \qquad
+ \int_0^t \textrm{ds}
\lla \fhat(s) \ot \ve(f), \widetilde{W}^*_s
\tilde{j}_s(\wt{\theta}(a)) \widetilde{W}_s (\ghat(s) \ot \ve(g)) \rra.
\end{align*}
where $\wt{\theta} = \big(\theta\ot\id_{\alg}\big)\circ\Com$.
This completes the proof.
\end{proof}

For each $t\geq 0$ denote the orthogonal projection from $\FFock$ onto
$\FFock_{[t, \infty[}$ by $P_{\kil,[t,\infty[}$. The following result
may be proved by differentiation, as with its predecessor for standard QS
cocycles, Lemma 4.2 of \cite{Stine}.

\begin{propn} \label{Stingen}
Let $\kil$ be an orthogonal direct sum of Hilbert spaces:
$\kil_0\oplus\kil_1$, let
$ \varphi\in CB\big(\alg;B(\wh{\kil}_0)\big)$ and
$\psi\in CB\big(\alg; B(\wh{\kil})\big)$, and
let  $k^0=l^{\varphi} \in
%\pPProcAtoCk0
\salvation$ and
$k=l^{\psi}\in \PProcAtoCk$
be the respective weak QS convolution cocycles. Then
\[k_t(a) = k^0_t(a) \odot P_{\kil_1,  [t,\infty[}
\;\;a\in \alg, t\geq 0,\]
if and only if
\[
\psi (a) = \left[\begin{matrix}\varphi(a) & 0 \\
 0 & - \Cou(a)I_1 \end{matrix} \right],
\;\;a\in \alg,
\]
where $I_1 = I_{\kil_1}$.
\end{propn}

We are ready for the main theorem of this section.

\begin{thm} \label{Stin}
Let $k \in     %\pPProcAtoCk0
\salvation$ be a Markov-regular CPC QS
convolution cocycle on a $C^*$-bialgebra $\alg$. Then there exists another
Hilbert space $\kil_1$, a Markov-regular, *-homomorphic
QS convolution cocycle $j \in \PProcAtoCk$,
where $\kil := \kil_0\oplus\kil_1$,
and a contractive process $W\in \ProcEk$, such that
\[ \wt{k}_t(a) = W_t^*j_t(a)W_t, \;\; t\geq 0, a \in \alg,\]
where $\wt{k}_t(a) := k_t(a) \ot P_{\kil_1,  [t, \infty[} $. A process $W$ may be chosen so that
it satisfies the QS differential equation
\begin{equation} \label{pert1} {\rm d}W_t = (l \otimes U_t W_t ){\rm d}
\Lambda_t,\;\;\; W_0 = I_{\FFock_\kil} \end{equation}
for some $l \in B(\hat{\kil})$ in which $U \in \ProcEk$
is the  projection-valued process given by
$U_t = j_t(1)$, $t \geq 0$.
\end{thm}

\begin{proof}
Let $\varphi\in CB\big(\alg; B(\kilhat)\big)$ be the stochastic  generator of $k$
(so that $k=l^{\varphi}$) and let
$(\Kil, \rho, D, \xi, d, e, t)$ be an associated tuple, as in
Theorem~\ref{spec}. Set $\kil_1=\Kil$ and define
$\theta: \alg \to B(\widehat{\kil})$ by
\[ \theta(a) = \left[ \begin{matrix}
\lambda(a) - t \Cou(a) & 0 & \delta^\dagger(a)\\
 0 & - \Cou(a) I_0 & 0 \\
\delta(a)  & 0 & \rho(a) - \Cou(a) I_1
\end{matrix}
 \right], \;\;\; a\in \alg, \]
where $I_i$ denotes $I_{\kil_i}$, $i=0,1$ and $\delta$ is the ($\rho,\Cou$)-derivation
$a\mapsto\big(\rho(a)-\Cou(a)I_{\kil}\big)|\xi\ra$.
The map $\theta$ is completely bounded and as such generates
a Markov-regular weak QS convolution cocycle
$j=l^{\theta} \in \PProcAtoCk$. It is easily checked that
$\theta$ satisfies the structure relations of Theorem~\ref{mult},
so $j$ is *-homomorphic.
Now choose any contraction $B\in B(\kil_1; \kil_0)$  and define
\[ \tau =   \left[ \begin{matrix} \frac{1}{2}t & \la \xi| & 0 \\
 0 & -I_0 & B  \\ 0 & D & -I_1     \end{matrix}\right]
\in B(\kilhat). \]
Then
\[ \tau^* + \tau + \tau^*\QSproj \tau = \left[ \begin{matrix} t & \la\xi| & 0 \\
 |\xi\ra & D^*D-I_0 & 0  \\ 0 & 0 & B^*B-I_1
\end{matrix}\right]\leq 0 ,\]
 as $B$ is  a contraction, and $\varphi(1)\leq 0$ (see (\ref{contrac})).

Theorem \ref{gener} yields the existence of  a contractive process $W \in \ProcEk$
satisfying the QS differential equation (\ref{pert1})
and  shows that the process $l \in \PProcAtoCk$ given by
\[ l_t(a) = W_t^*j_t(a)W_t, \;\; t\geq 0, a \in \alg,\]
is equal to $l^{\psi}$ where  $\psi:\alg \to B(\hat{\kil})$ is defined by
\begin{eqnarray*}
\psi(a) &=& \Cou(a) \left( \tau^*+\tau + \tau^*\QSproj \tau \right) + (1 + \tau^* \QSproj ) \theta(a) (1 + \QSproj \tau) \\
 &=& \Cou(a) \left[ \begin{matrix} t & \la\xi| & 0 \\
 |\xi\ra & D^*D-I_0 & 0  \\ 0 & 0 & B^*B-I_1     \end{matrix}\right] \\
  &&+ \left[ \begin{matrix} 1 & 0 & 0 \\
 0 & 0& D^* \\ 0 & B^* & 0     \end{matrix}\right] \cdot  \left[ \begin{matrix}   \lambda(a) - t \Cou(a) & 0 &
\delta^\dagger(a)\\
 0 & - \Cou(a) I_0 & 0 \\ \delta(a) & 0 & \rho(a) - \Cou(a) I_1
\end{matrix}\right]  \cdot  \left[ \begin{matrix} 1 & 0 & 0 \\
 0 &  0 & B \\ 0 & D & 0   \end{matrix}\right] \\
 &=&  \Cou(a) \left[ \begin{matrix} t & \la\xi| & 0 \\
 |\xi\ra & D^*D-I_0 & 0  \\ 0 & 0 & B^*B-I_1     \end{matrix}\right] \\
 && + \left[ \begin{matrix}   \lambda(a) - t \Cou(a) & \delta^\dagger(a)D & 0 \\
 D^* \delta(a) & D^* \rho(a) D - \Cou(a) D^*D & 0 \\ 0 & 0 &  - \Cou(a) B^*B \end{matrix} \right] \\
&=& \left[ \begin{matrix}   \lambda(a) &  \delta^\dagger(a)D + \Cou(a)\la \xi | & 0 \\
 D^* \delta(a)+ \Cou(a)|\xi  \ra & D^* \rho(a) D - \Cou(a)I_0 & 0 \\ 0 & 0
&  - \Cou(a) I_1 \end{matrix}\right] =
\left[ \begin{matrix}   \varphi(a) & 0 \\  0 &  - \Cou(a) I_{\kil}
\end{matrix}\right]. \end{eqnarray*}
Application of Proposition~\ref{Stingen} now completes the proof.
\end{proof}

\section*{Acknowledgements}
The author would like to express his gratitude to Martin Lindsay for suggesting the problem
considered in this paper and for generous comments improving its final form.
 The work was partially supported by
the Polish KBN Research Grant 2P03A 03024.


\begin{thebibliography}{GLSW}



\bibitem[Acc] {Acc}
L.\,Accardi,
On the quantum Feynman-Kac formula,
\emph{Rend. \!Sem. \!Mat. \!Fis. \!Milano}
\textbf{48} (1978), 135--180.

\bibitem[Bel] {Slava}
V.P.\,Belavkin, Quantum stochastic positive evolutions: characterization,
construction, dilation,
\emph{Comm.\,Math.\,Phys.}  \textbf{184}  (1997)  no.\,3, 533--566.

\bibitem[BlH] {hyp}
W.R.\,Bloom, and H.\,Heyer, ``Harmonic analysis of probability measures on hypergroups,''
de Gruyter Studies in Mathematics, 20. Walter de Gruyter \& Co., Berlin, 1995.


\bibitem[EfR] {ERuan}
E.G.\,Effros  and Z.J.\,Ruan, ``Operator Spaces,'' Oxford University Press,
Oxford 2000.

\bibitem[ChV]{ChaV}
Yu.\,Chapovsky and L.\,Vainerman,
Compact quantum hypergroups,
\emph{J.\,Operator Theory} \textbf{41} (1999) no.\,2, 261--289.

\bibitem[Fra]{franz}
U.\,Franz,
L\'evy processes on quantum groups and dual groups,
\emph{in} ``Quantum Independent Increment Processes,
 Vol. \!II: Structure of Quantum L\'evy Processes,
Classical Probability and Physics,''
\emph{eds.\ U. \!Franz \& M. \!Sch\"urmann},
Lecture Notes in Mathematics \textbf{1866},
Springer-Verlag, Heidelberg 2005.


\bibitem[FrS] {UweSch}
U.\,Franz and M.\,Sch\"{u}rmann,
L\'evy processes on quantum
hypergroups, \emph{in}
``Infinite Dimensional
Harmonic Analysis,''
\emph{eds.\ H. \!Heyer, T. \!Hirai \& N. \!Obata}
Gr\"abner, Altendorff 2000,
pp. \!93--114.


\bibitem[GLW] {Stine}
D.\,Goswami, J.M.\,Lindsay and S.J.\,Wills, A stochastic Stinespring theorem,
\emph{Math.\,Ann.}
\textbf{319} (2001) no.\,4,  647--673.

\bibitem[GLSW] {dilate}
D.\,Goswami, J.M.\,Lindsay, K.B.\,Sinha  and S.J.\,Wills, Dilation of Markovian cocycles
on a von Neumann algebra, \emph{Pacific J.\,Math.}
\textbf{211} (2003), 221--247.

\bibitem[Kal] {Kal1}
A.A.\,Kalyuzhnyi,
Conditional expectations on compact quantum groups and new examples of quantum hypergroups, \emph{Methods Funct.\,Anal.\,Topology}
\textbf{7}  (2001)  no.\,4, 49--68.


\bibitem[KaC] {Kal2}
A.A.\,Kalyuzhnyi and Yu.A.\,Chapovsky,
A factorization of conditional expectations on Kac algebras and quantum double coset hypergroups,
\emph{Ukrain.\,Mat.\,Zh.}  \textbf{55}  (2003)  no.\,12, 1669--1677;  translation in  \emph{Ukrainian Math.\, J.} \textbf{55}
(2003)  no. 12, 1994--2005.

\bibitem[Lin]{lect}
J.M.\,Lindsay,
Quantum stochastic analysis -- an introduction,
\emph{in} ``Quantum Independent Increment Processes,
Vol. \!I: From Classical Probability to Quantum Stochastics,"
\emph{eds.\ U. \!Franz \& M. \!Sch\"urmann},
Lecture Notes in Mathematics \textbf{1865},
Springer-Verlag, Heidelberg 2005.

\bibitem[LiP]{lp}
J.M.\,Lindsay and K.R.\,Parthasarathy,
On the generators of quantum stochastic flows,
\emph{J.\,Funct.\,Anal.}
\textbf{158} (1998), 521--549.


\bibitem[$\text{LS}_1$]
{LSqscc1}
J.M.\,Lindsay and A.G.\,Skalski,
Quantum stochastic convolution cocycles I,
\emph{Ann. \!Inst. \!H. \!Poincar\'{e}, Probab. \!Statist.}
\textbf{41} (2005) no.\,3  (En hommage \`a Paul-Andr\'e Meyer),
581--604.

\bibitem[$\text{LS}_2$]
%\bibitem[LS]
{LSbedlewo}
J.M.\,Lindsay and A.G.\,Skalski,
Quantum stochastic convolution cocycles --- algebraic and $C^*$-algebraic,
\emph{Banach Center Publ.} \textbf{73} (2006), 313--324.


\bibitem[$\text{LS}_3$]
%\bibitem[LS]
{LSqsde}
J.M.\,Lindsay and A.G.\,Skalski,
On quantum stochastic differential equations, \emph{J.\,Math.\,Anal.\,Appl.} (2007), doi:10.1016/j.jmaa.2006.07.105

\bibitem[$\text{LS}_4$]
%\bibitem[LS]
{LSqscc2}
J.M.\,Lindsay and A.G.\,Skalski,
Quantum stochastic convolution cocycles II, \emph{preprint}.

\bibitem[LW$_1$]{lw1}
J.M.\,Lindsay and S.J.\,Wills,
Existence, positivity, and contractivity for quantum stochastic flows
with infinite dimensional noise,
\emph{Probab.\,Theory Related Fields}
\textbf{116} (2000), 505--543.

\bibitem[LW$_2$]{lwjfa}
J.M.\,Lindsay and S.J.\,Wills,
Markovian cocycles on operator algebras, adapted to a Fock filtration,
\emph{J.\,Funct.\,Anal.}
\textbf{178} (2000) no.\,2, 269--305.

\bibitem[LW$_3$]{lwblms}
J.M.\,Lindsay and S.J.\,Wills,
Existence of Feller cocycles on a $C^*$-algebra,
\emph{Bull.\,London Math.\,Soc.}
\textbf{33} (2001) no.\,5, 613--621.

\bibitem[LW$_4$]{lwjlms}
J.M.\,Lindsay and S.J.\,Wills,
Homomorphic Feller cocycles on a $C^*$-algebra,
\emph{J. \!London Math. \!Soc.} (\emph{2})
\textbf{68} (2003) no. \!1, 255--272.

\bibitem[Sch]{schu}
M.\,Sch\"{u}rmann,
``White Noise on Bialgebras,''
Lecture Notes in Mathematics \textbf{1544},
Springer, Heidelberg 1993.

\end{thebibliography}
\end{document}